\input amstex
\input amsppt.sty
\magnification=\magstep1
\hsize=30truecc
\vsize=22.2truecm
\baselineskip=16truept
\nologo
\TagsOnRight
\def\N{\Bbb N}
\def\Z{\Bbb Z}

\def\Q{\Bbb Q}

\def\C{\Bbb C}
\def\l{\left}
\def\r{\right}
\def\bg{\bigg}
\def\({\bg(}
\def\[{\bg[}
\def\){\bg)}
\def\]{\bg]}
\def\t{\text}
\def\f{\frac}

\def\mo{\roman{mod}}
\def\per{\roman{per}}

\def\em{\emptyset}
\def\se {\subseteq}
\def\sp {\supseteq}
\def\sm{\setminus}

\def\bi{\binom}
\def\eq{\equiv}
\def\cs{\cdots}
\def\ls{\leqslant}
\def\gs{\geqslant}
\def\al{\alpha}

\def\da{\delta}

\def\bi{\binom}
\def\Proof{\noindent{\it Proof}}

\def\Remark{\medskip\noindent{\it  Remark}}

 \topmatter
 \hbox{Preprint}
 \medskip
 \title Connections between covers of $\Z$ and subset sums\endtitle
 \rightheadtext{Connections between covers of $\Z$ and subset sums}
 \author Zhi-Wei Sun\endauthor
 \affil Department of Mathematics, Nanjing University
 \\Nanjing 210093, People's Republic of China
 \\  zwsun\@nju.edu.cn
 \\ {\tt http://maths.nju.edu.cn/$\sim$zwsun}\endaffil
 \medskip

\abstract In this paper we establish connections between covers of $\Z$ by residue classes and subset sums in a field. Suppose that $\{a_s(n_s)\}_{s=0}^k$ covers each integer at least $p$ times
with the residue class $a_0(n_0)$ irredundant, where $p$ is a prime
not dividing any of $n_1,\ldots,n_k$.
Let $m_1,\ldots,m_k\in\Z$ be relatively prime to
$n_1,\ldots,n_k$ respectively. For any
$c,c_1,\ldots,c_k\in\Z/p\Z$ with $c_1\cdots c_k\not=0$, we show that the set
$$\bigg\{\bg\{\sum_{s\in I}\f{m_s}{n_s}\bg\}:\, I\subseteq\{1,\ldots,k\}
\ \t{and}\ \sum_{s\in I}c_s=c\bigg\}$$
contains an arithmetic progression of length $n_0$ with common
difference $1/n_0$, where $\{x\}$ denotes the fractional part of a real number $x$.

\endabstract
\thanks 2020 {\it Mathematics Subject Classification}.
Primary 11B25, 11B75; Secondary 05E99, 11B13, 11C08, 11D68, 11T99.
\newline\indent Supported by
the National Natural Science Foundation of China (grant 11971222).
\endthanks
\endtopmatter
 \document
\heading {1. Introduction} \endheading

 For a finite set $S=\{a_1,\ldots,a_k\}$ contained in the ring $\Z$ or a field,
 sums in the form $\sum_{s\in I}a_s$ with $I\se[1,k]=\{1,\ldots,k\}$ are
 called {\it subset sums} of $S$. It is interesting to provide a lower bound for
 the cardinality of the set
 $$\{a_1x_1+\cs+a_kx_k:\, x_1,\ldots,x_k\in\{0,1\}\}
 =\bg\{\sum_{s\in I}a_s:\, I\se[1,k]\bg\}.$$
 A more general problem is to study restricted sumsets
 in the form
 $$\{x_1+\cs+x_k:\, x_1\in X_1,\ldots, x_k\in X_k,\
 P(x_1,\ldots,x_k)\not=0\}\tag1.1$$
 where $X_1,\ldots,X_k$ are subsets of a field
 and $P(x_1,\ldots,x_k)$ is a polynomial with coefficients in the
 field.

 Let $p$ be a prime. In 1964 Erd\H os and Heilbronn [EH] conjectured that
 if $\em\not=X\se\Z_p=\Z/p\Z$ then
 $$|\{x_1+x_2:\, x_1,x_2\in X\ \t{and}\ x_1\not=x_2\}|\gs\min\{p,\,2|X|-3\}.$$
 This conjecture was first confirmed by Dias da Silva and
 Hamidoune [DH] in 1994, who obtained
 a generalization which implies that
 if $S\se\Z_p$ and $|S|>\sqrt{4p-7}$ then
 any element of $\Z_p$ is a subset sum of $S$.
 In this direction the most powerful tool is
 the following remarkable principle (see Alon [A99, A03]) rooted in
 Alon and Tarsi [AT] and  applied in
 [AF], [ANR1, ANR2], [DKSS], [HS], [LS], [PS], [S03b], [S08] and [SZ].

 \proclaim{Combinatorial Nullstellensatz {\rm (Alon [A99])}}
 Let $X_1,\ldots,X_k$ be finite subsets
 of a field $F$ with $|X_s|>l_s$ for $s\in[1,k]$
 where $l_1,\ldots,l_k\in\N=\{0,1,\ldots\}$.
 If $f(x_1,\ldots,x_k)\in F[x_1,\ldots,x_k]$, and
  $[x_1^{l_1}\cs x_k^{l_k}]f(x_1,\ldots,x_k)$
 $($the coefficient
 of the monomial $\prod_{s=1}^kx_s^{l_s}$ in $f)$
 is nonzero and $\sum_{s=1}^kl_s$ is
 the total degree of $f$,
 then there are $x_1\in X_1,\ldots,x_k\in X_k$ such that
 $f(x_1,\ldots,x_k)\not=0$.
 \endproclaim

 One of many applications of the Combinatorial Nullstellensatz
 is the following result of [AT] concerning a conjecture of
 J\" ager.

 \proclaim{Alon-Tarsi Theorem} Let $F$ be a finite field with
 $|F|$ not a prime, and let $M$ be a nonsingular $k\times k$
 matrix over $F$. Then there exists a vector $\vec x=(x_1,\ldots,x_k)^T$
 with $x_1,\ldots,x_k\in F$
 such that neither $\vec x$ nor $M\vec x$ has zero component.
 \endproclaim

 Now we turn to covers of $\Z$ by finitely many residue classes.

 For $a\in\Z$ and $n\in\Z^+=\{1,2,3,\ldots\}$, set
 $$a(n)=a+n\Z=\{a+nx:\, x\in\Z\}$$
 and call it a residue class with modulus $n$. For a finite system
 $$A=\{a_s(n_s)\}_{s=1}^k\tag1.2$$
 of residue classes, we define its {\it covering function} $w_A:\Z\to\Z$ by
 $$w_A(x)=|\{1\ls s\ls k:\, x\in a_s(n_s)\}|.$$
 For properties of the covering function $w_A(x)$, one can consult [S03a, S04].
 As in Sun [S97, S99], we
 call $m(A)=\min_{x\in\Z}w_A(x)$
 the {\it covering multiplicity} of (1.2).

 Erd\H os [E50] initiated the study of covers of $\Z$ by residue classes.
 Zhang [Z89] showed that if $(1.2)$ covers all the integers then $\sum_{s\in I}1/n_s\in\Z^+$
 for some $I\se\{1,\ldots,k\}$.
 Let $m\in\Z^+$. We call (1.2) an {\it $m$-cover}
 of $\Z$ if $m(A)\ge m$.
 If (1.2) forms an $m$-cover of $\Z$
 but $A_t=\{a_s(n_s)\}_{s\not=t}$ does not,
 then we say that $(1.2)$ is an $m$-cover of $\Z$
 with $a_t(n_t)$ {\it essential}. 
 
 The author [S99] established the following result.
 \proclaim{Theorem {\rm (Sun [S99])}} Let $m\in\Z^+$ and let $A=\{a_s(n_s)\}_{s=0}^k$
 be an $m$-cover of $\Z$ with $a_0(n_0)$ essential. Let $m_1,\ldots,m_k\in\Z^+$
 be relatively prime to $n_1,\ldots,n_s$ respectively. Then the set
 $$\bg\{\bg\{\sum_{s\in I}\f{m_s}{n_s}\bg\}:\, I\se[1,k]\
\t{and}\ \bg\lfloor\sum_{s\in I}\f{m_s}{n_s}\bg\rfloor\gs m-1\bg\}$$
contains an arithmetic progression of length $n_0$
with common difference $1/n_0$, where $\{x\}$ and $\lfloor x\rfloor$
are the fractional part and the integer part of a real number $x$.
\endproclaim

 Subset sums seem to have nothing to do with
 covers of $\Z$. Before our work no one else has realized their close
 connections. Can you imagine that the Alon-Tarsi theorem
 are related to covers of $\Z$?

    The purpose of this paper is to present a surprising unified approach and
 embed the study of subset sums
 in the investigation of covers. The key point of our unification
 is to compare the following two sorts of quantities:

(a) Degrees of multi-variable polynomials over rings or fields,

(b) Covering multiplicities of covers of $\Z$ by residue classes.

 In Section 2 we will present a general unified theorem
 connecting subset sums with covers of $\Z$
 and derive from it some consequences.

 In Section 3 we will pose a formula for polynomials over a ring.

 On the basis of Section 3, the reader will understand quite well
 the technique in Section 4 used to prove Theorem 2.1
 which connects covers of $\Z$ with subset sums.
 
 The author [S03c] announced a unified approach to covers of $\Z$, subset sums and zero-sum problems.
The detailed connections between covers of $\Z$ and zero-sum problems were published in [S09]. 
 
Let $m\in\Z^+$. The system (1.2) is called an {\it $m$-system} if $w_A(x)\ls m$ for all $x\in\Z$.
One may wonder whether such systems are
also related to subset sums.
Let
$$A^*=\{a_s+r(n_s):\, r\in[1,n_s-1]\ \t{and}\ s\in[1,k]\}\tag1.3$$
and call it {\it the dual system} of (1.2) as in [S10].
Then $w_A(x)+w_{A^*}(x)=k$ for all $x\in\Z$. Thus (1.2) is an $m$-system
if and only if $m(A^*)\gs k-m$.
In light of this, we can reformulate our results related to covers of $\Z$
in terms of $m$-systems.

\heading{2. A General Theorem and its Consequences}\endheading

Now we state our general theorem connecting covers of $\Z$ with subset sums.

 \proclaim{Theorem 2.1} Let $A_0=\{a_s(n_s)\}_{s=0}^k$ be a system of residue classes
 with $w_{A_0}(a_0)=m(A_0)$.
 Let $m_1,\ldots,m_k\in\Z$ be relatively prime to $n_1,\ldots,n_k$ respectively.
 Let $J\se\{1\ls s\ls k:\, a_0\in a_s(n_s)\}$ and
 $P(x_1,\ldots,x_k)\in F[x_1,\ldots,x_k]$
 where $F$ is a field with characteristic not dividing $N=[n_1,\ldots,n_k]$.
 Assume that $0\ls\deg P\ls|J|$ and
 $$\[\prod_{j\in J}x_j\]P(x_1,\ldots,x_k)
 (x_1+\cs+x_k)^{|J|-\deg P}\not=0.\tag2.1$$
 Let $X_1=\{b_1,c_1\},\ldots,X_k=\{b_k,c_k\}$
 be subsets of $F$ such that $b_s=c_s$ only if $a_0\in a_s(n_s)$ and $s\not\in J$.
 Then, for some $0\ls\alpha<1$, we have
 $$|S_r|\gs |J|-\deg P+1>0\quad\t{for all}\ r=0,1,\ldots,n_0-1,\tag2.2$$
 where $$S_r=\bg\{\sum_{s=1}^kx_s:\, x_s\in X_s,\ P(x_1,\ldots,x_k)\not=0,
 \  \bg\{\sum^k\Sb s=1\\x_s\not=b_s\endSb\f{m_s}{n_s}\bg\}=\f {\al+r}{n_0}\bg\}.
 \tag2.3$$
 \endproclaim

  In the case $n_0=n_1=\cs=n_k=1$, Theorem 2.1
 yields the following basic lemma of the so-called polynomial
 method due to Alon, Nathanson and Ruzsa [ANR1, ANR2]:
 Let $X_1,\ldots,X_k$ be subsets of a field $F$ with
 $|X_s|>l_s\in\{0,1\}$ for $s\in[1,k]$. If $P(x_1,\ldots,x_k)
 \in F[x_1,\ldots,x_k]\sm\{0\}$, $\deg P\ls \sum_{s=1}^kl_s$ and
 $$[x_1^{l_1}\cdots
 x_k^{l_k}]P(x_1,\ldots,x_k)(x_1+\cs+x_k)^{\sum_{s=1}^kl_s-\deg P}
 \not=0,$$
 then
 $$\bg|\bg\{\sum_{s=1}^kx_s:\, x_s\in X_s\ \t{and}
 \ P(x_1,\ldots,x_k)\not=0\bg\}\bg|
 \gs\sum_{s=1}^kl_s-\deg P+1.$$
 Actually this remains valid even if $l_s$ may be greater than one.

\proclaim{Corollary 2.1} Let $A_0=\{a_s(n_s)\}_{s=0}^k$ be an $m$-cover of
$\Z$ with $a_0(n_0)$ essential. Let $m_1,\ldots,m_k\in\Z$
be relatively prime to $n_1,\ldots,n_k$ respectively.
Let $F$ be a field with characteristic $p$ not dividing $[n_1,\ldots,n_k]$, and let
$X_1=\{b_1,c_1\},\ldots,X_k=\{b_k,c_k\}$
be any subsets of $F$ with cardinality $2$. Then, for some
$0\ls\al<1$, we have
$$\bg|\bg\{\sum_{s=1}^{k}x_s:\, x_s\in X_s,
\ \bg\{\sum\Sb 1\ls s\ls k\\x_s=c_s\endSb\f{m_s}{n_s}\bg\}=\f{\al+r}{n_0}\bg\}\bg|
\gs\min\{p',m\}\tag2.4$$
for all $r\in[0,n_0-1]$,
where $p'=p$ if $p$ is a prime, and $p'=+\infty$ if $p=0$.
\endproclaim
\Proof. Since $a_0(n_0)$ is essential, there is an $a\in a_0(n_0)$ such that
$w_{A_0}(a)=m$. Note that $a_0(n_0)=a(n_0)$. Choose
$J\se\{1\ls s\ls k:\, a\in a_s(n_s)\}$ with
$|J|=\min\{p',m\}-1$. Then
$$\[\prod_{j\in J}x_j\](x_1+\cs+x_k)^{|J|}
=\[\prod_{j\in J}x_j\]|J|!\prod_{j\in J}x_j\not=0$$
since $|J|<p'$. Now it suffices to apply
Theorem 2.1 with $P(x_1,\ldots,x_k)=1$. \qed

\Remark\ 2.1. Let $A_0=\{a_s(n_s)\}_{s=0}^k$ be an $m$-cover of $\Z$
with $a_0(n_0)$ essential. And let $m_1,\ldots,m_k\in\Z^+$ be relatively prime to
$n_1,\ldots,n_k$ respectively.
By Corollary 2.1 in the case $F=\Q$ and $X_s=\{0,m_s/n_s\}\ (1\ls s\ls k)$,
for some $0\ls\al<1$ we have
$$\min_{r\in[0,\, n_0-1]}\bg|\bg\{\sum_{s\in I}\f{m_s}{n_s}:
\, I\se[1,k]\ \t{and}\ \bg\{\sum_{s\in I}\f{m_s}{n_s}\bg\}
=\f{\al+r}{n_0}\bg\}\bg|\gs m.\tag2.5$$
This implies that
$$\bg\{\bg\{\sum_{s\in I}\f{m_s}{n_s}\bg\}:\, I\se[1,k]\
\t{and}\ \bg\lfloor\sum_{s\in I}\f{m_s}{n_s}\bg\rfloor\gs m-1\bg\}$$
contains an arithmetic progression of length $n_0$
with common difference $1/n_0$,
which was first established by the author in [S99].
In 2007 the author [S07] showed that
if the covering function $w_{A_0}(x)$ is periodic modulo $n_0$ then
(2.5) holds with $m_1=\cdots=m_k=1$ and $\al=0$.
\medskip

Inspired by Corollary 2.2 (first announced in [S03c])
and an earlier paper [S97], the author [S10] proved
prove that if $A_0=\{a_s(n_s)\}_{s=0}^k$ forms an $m$-cover of $\Z$ with
$\sum_{s=1}^k1/n_s<m$ then for any $a=0,1,2,\ldots$ we have
$$\bg|\bg\{I\se[1,k]:\, \sum_{s\in I}\f1{n_s}=\f a{n_0}\bg\}\bg|\gs\bi{m-1}{\lfloor a/n_0\rfloor}.$$

\proclaim{Corollary 2.2}
Let $A_0=\{a_s(n_s)\}_{s=0}^k$ be a $p$-cover of $\Z$
with $a_0(n_0)$ essential, where $p$ is a prime
not dividing any of $n_1,\ldots,n_k$.
Let $m_1,\ldots,m_k\in\Z$ be relatively prime to
$n_1,\ldots,n_k$ respectively. Then, for any
$c,c_1,\ldots,c_k\in\Z_p$ with $c_1\cdots c_k\not=0$, the set
$$\bg\{\bg\{\sum_{s\in I}\f{m_s}{n_s}\bg\}:\, I\se[1,k]
\ \t{and}\ \sum_{s\in I}c_s=c\bg\}\tag2.6$$
contains an arithmetic progression of length $n_0$ with common
difference $1/n_0$.
\endproclaim
\Proof. By Corollary 2.1 in the case $F=\Z_p=\Z/p\Z$ and $X_s=\{0,c_s\}\ (1\ls s\ls k)$,
for some $0\ls\al<1$ we have
$\{\sum_{s\in I}c_s:\,\{\sum_{s\in I}m_s/n_s\}=(\al+r)/n_0\}=\Z_p$
for every $r\in[0,n_0-1]$. So the desired result follows.
\qed

\Remark\ 2.2. The author's colleague Z. Y. Wu
 once asked whether for any prime $p$ and  $c_1,\ldots,c_{p-1}
 \in\Z_p\sm\{0\}$ there is an $I\se[1,p-1]$ such that
 $\sum_{s\in I}c_s=1$.
 Corollary 2.2 in the case $n_0=n_1=\cs=n_k=1$,
 provides an affirmative answer to this question.

\proclaim{Corollary 2.3} Let $A_0=\{a_s(n_s)\}_{s=0}^k$ be an $m+1$-cover of $\Z$
with $m\in\N$ and $w_{A_0}(a_0)=m+1$.
Let $F$ be a field with characteristic not dividing any of $n_1,\ldots,n_k$,
and let $X_1,\ldots,X_k$ be subsets of $F$
with cardinality $2$. Let $a_{ij},b_i\in F$ and $c_i\in X_i$
for all $i\in[1,m]$ and $j\in[1,k]$.
If $m_1,\ldots,m_k\in\Z$ are relatively prime to $n_1,\ldots,n_k$
respectively, and
$$\per(a_{ij})_{i\in[1,m],\, j\in J}:=\sum_{\{j_1,\ldots,j_m\}=J}
a_{1j_1}\cdots a_{mj_m}\not=0$$
where $J=\{1\ls s\ls k:\, a_0\in a_s(n_s)\}$,
then the set
$$\bg\{\bg\{\sum\Sb 1\ls s\ls k\\x_s=c_s\endSb\f{m_s}{n_s}\bg\}:\,
x_s\in X_s\ \t{and}\ \sum_{j=1}^ka_{ij}x_j\not=b_i\ \t{for all}
\ i\in[1,m]\bg\}\tag2.7$$
contains an arithmetic progression of length $n_0$
with common difference $1/n_0$.
\endproclaim
\Proof. Note that $|J|=m$. Set
$P(x_1,\ldots,x_k)=\prod_{i=1}^m(\sum_{j=1}^ka_{ij}x_j-b_i)$.
Then
$$\[\prod_{j\in J}x_j\]P(x_1,\ldots,x_k)=\[\prod_{j\in J}x_j\]
\prod_{i=1}^m\sum_{j\in J}a_{ij}x_j=\per(a_{ij})_{i\in[1,m],\,j\in J}\not=0.$$
In view of Theorem 2.2, the set
$$\bg\{\bg\{\sum\Sb 1\ls s\ls k\\x_s=c_s\endSb\f{m_s}{n_s}\bg\}:\, x_s\in X_s,
\ P(x_1,\ldots,x_k)\not=0\bg\}$$
contains $\{(\al+r)/n_0:\,r\in[0,n_0-1]\}$ for some $0\ls\al<1$.
We are done. \qed

\Remark\ 2.3. When $n_0=n_1=\cs=n_k=1$,
Corollary 2.3 yields the useful permanent lemma of Alon [A99].

\proclaim{Corollary 2.4} Let $A_0=\{a_s(n_s)\}_{s=0}^k$
be an $m+1$-cover of $\Z$ with $a_0(n_0)$ essential.
Let $m_1,\ldots,m_k$ be integers relatively prime to
$n_1,\ldots,n_k$ respectively. Let $F$ be a field
of prime characteristic $p$, and let $a_{ij},b_i\in F$ for
all $i\in[1,m]$ and $j\in[1,k]$.
Set
$$X=\bg\{\sum_{j=1}^kx_j:\,
x_j\in[0,p-1]\ \t{and}\ \sum_{j=1}^kx_ja_{ij}\not=b_i
\ \t{for all}\ i\in[1,m]\bg\}.\tag2.8$$
If $p$ does not divide $N=[n_1,\ldots,n_k]$
and the matrix $(a_{ij})_{1\ls i\ls m,\,1\ls j\ls k}$
has rank $m$, then the set
$$S:=\bg\{\bg\{\sum_{s\in I}\f{m_s}{n_s}\bg\}:\, I\se[1,k]\ \t{and}\
|I|\in X\bg\}\tag2.9$$
contains an arithmetic progression of length $n_0$
with common difference $1/n_0$;
in particular, when $n_0=N$ we have
$S=\{r/N:\, r\in[0,N-1]\}$.
\endproclaim
\Proof. As $a_0(n_0)$ is essential,
for some $a\in a_0(n_0)$ we have $w_{A_0}(a)=m+1$.
Without loss of generality we assume that
the matrix $M=(a_{ij})_{i,j\in[1,m]}$
is nonsingular, and that
$\{1\ls s\ls k:\, a\in a_s(n_s)\}=[1,m]$ (otherwise
we can rearrange the $k$ residue classes in a suitable order).

Since $\det M\not=0$,
by [AT] there are $l_1,\ldots,l_m
\in[0,p-1]$ with $l_1+\cs+l_m=m$ such that
$\per(M^*)\not=0$, where $M^*$ is an $m\times m$ matrix
whose columns consist of $l_1$ copies of the first column of
$M$, $\ldots$, $l_m$ copies of the $m$th column of $M$.
Let $e$ denote the identity of the field $F$.
By Corollary 2.5, there exists $0\ls\al<1$ such that
for any $r\in[0,n_0-1]$ there are $\da_1,\ldots,\da_k\in\{0,e\}$
for which $\{\sum_{\da_s=e}m_s/n_s\}=(\al+r)/n_0$
and
$$\sum_{j=1}^ka_{ij}(\da_{l_1+\cs+l_{j-1}+1}
+\cs+\da_{l_1+\cs+l_j})\not=b_i\quad \t{for all}\ i\in[1,m],$$
where $l_j=1$ for any $j\in[m+1,k]$.
Observe that
$$x_j=|\{l_1+\cs+l_{j-1}<s\ls l_1+\cs+l_j:\,
\da_s=e\}|\ls l_j<p$$
and $|\{1\ls s\ls k:\,\da_s=e\}|=x_1+\cs+x_k\in X$.
So the set $S$ given by (2.13) contains $\{(\al+r)/n_0:\, r\in[0,n_0-1]\}$.
As $T=\{r/N:\, r\in[0,N-1]\}\sp S$,
we have $S=T$ if $n_0=N$.
This concludes the proof.
\qed

\Remark\ 2.4. The Alon-Tarsi Theorem stated in Section 1
follows from Corollary 2.4
for the following reason:
Let $F$ be a field of prime characteristic $p$ with identity $e$.
If the matrix $(a_{ij})_{i,j\in[1,k]}$ over $F$ is nonsingular
and $c\in F\sm\{0,e,\ldots,(p-1)e\}$, then by Corollary 2.4 in the
case $n_0=n_1=\cs=n_k=1$ there are $x_1,\ldots,x_k\in[0,p-1]$
with $x_1+\cs+x_k\ls k$ such that
$$\sum_{j=1}^kx_ja_{ij}\not=-c\sum_{j=1}^ka_{ij},
\quad\ \t{i.e.}\ \sum_{j=1}^ka_{ij}x_j^*\not=0$$
for all $i\in[1,k]$ where $x_j^*=x_je+c\not=0$.
\medskip

\proclaim{Corollary 2.5} Let $(1.2)$ be an $m$-cover of $\Z$
with $a_k(n_k)$ essential and $n_k=N_A$. Let $m_1,\ldots,m_{k-1}
\in\Z$ be relatively prime to $n_1,\ldots,n_{k-1}$ respectively.
Then, for any $J\se K=\{1\ls s<k:\, a_k\in a_s(n_s)\}$ and $r\in[0,N_A-1]$,
there exists an $I\se[1,k-1]$ such that $I\cap K=J$
and $\{\sum_{s\in I}m_s/n_s\}=r/N_A$.
\endproclaim

\Proof. Clearly $w_A(a_k)=m$.
Fix $J\se K$. By Theorem 2.1 in the case $F=\Q$, the set
$$\align S=&\bg\{\bg\{\sum\Sb 1\ls s<k\\x_s\not=0\endSb\f{m_s}{n_s}\bg\}
:\, x_s\in\{0,1\},\ \prod_{j\in J}x_j\not=0,
\ x_s\in\{0\}\ \t{for}\ s\in K\sm J\bg\}
\\=&\bg\{\bg\{\sum_{s\in I}\f{m_s}{n_s}\bg\}:\,
I\se[1,k-1]\ \t{and}\ I\cap K=J\bg\}
\endalign$$
contains an arithmetic progression of length $n_k=N_A$
with common difference $1/n_k=1/N_A$.
Since $T=\{a/N_A:\, a\in[0,N_A-1]\}\sp S$,
we must have $S=T$ and thus
 $\{\sum_{s\in I}m_s/n_s\}=r/N_A$ for some $I\se[1,k-1]$
with $I\cap K=J$. \qed

\Remark\ 2.5. On the basis of the author's work [S95],
his brother Z.-H. Sun pointed out that if
(1.2) forms a cover of $\Z$ with $a_k(n_k)$
essential and $n_k=N_A$,
then $\{\{\sum_{s\in I}1/n_s\}:\,
I\se[1,k-1]\}=\{r/N_A:\, r\in[0,N_A-1]\}$.
This follows from Corollary 2.5 in the special case
$m=m_1=\cdots=m_{k-1}=1$.
\medskip

 Let $n>1$ be an integer, and let $m_1,\ldots,m_{n-1}\in\Z$ be
relatively prime to $n$.
Applying Corollary 2.5 to the trivial cover $\{r(n)\}_{r=0}^{n-1}$,
we find that the set $\{\sum_{s\in I}m_s:\, I\se[1,n-1]\}$
contains a complete system of residues modulo $n$.
This is more general than the positive answer
to Wu's question mentioned in Remark 2.2.

 \heading{3. A Useful Polynomial Formula and its Applications}\endheading
 
 For a predicate $P$, we define 
 $$\[\![P]\!]=\cases1&\t{if}\ P\ \t{holds},
 \\0&\t{otherwise}.\endcases$$
 
 The author [S03] first announed the following result in 2003. 

\proclaim{Theorem 3.1} Let $R$ be a  ring with identity, and let
 $f(x_1,\ldots,x_k)$ be a polynomial over $R$. If $J\se[1,k]$
and $|J|\gs\deg f$, then we have the formula
$$\sum_{I\se J}(-1)^{|J|-|I|}f([\![1\in I]\!],\ldots,[\![k\in I]\!])
=\[\prod_{j\in J}x_j\]f(x_1,\ldots,x_k).\tag3.1$$
\endproclaim

\Proof. Write
$f(x_1,\ldots,x_k)=\sum_{j_1,\ldots,j_k\gs0}c_{j_1,\ldots,j_k}\prod_{s=1}^kx_s^{j_s},$
and observe that if
$\em\not=J'\se[1,k]$ then $0=\prod_{j\in J'}(1-1)=\sum_{I\se J'}(-1)^{|I|}$.
Therefore
$$\align &\sum_{I\se J}(-1)^{|I|}f([\![1\in I]\!],\ldots,[\![k\in I]\!])
\\=&\sum_{I\se J}(-1)^{|I|}
\sum\Sb j_1,\ldots,j_k\gs0\\\{s:\, j_s\not=0\}\se I\endSb
c_{j_1,\ldots,j_k}
\\=&\sum\Sb j_1,\ldots,j_k\gs0\\\{s:\, j_s\not=0\}\se J\endSb
\sum_{\{s:\, j_s\not=0\}\se I\se J}(-1)^{|I|}c_{j_1,\ldots,j_k}
\\=&\sum\Sb j_1,\ldots,j_k\gs0\\\{s:\, j_s\not=0\}\se J\endSb
\sum_{I'\se J\sm\{s:\, j_s\not=0\}}
(-1)^{|I'|}(-1)^{|\{s:\, j_s\not=0\}|}c_{j_1,\ldots,j_k}
\\=&\sum\Sb j_1,\ldots,j_k\gs0\\\{s:\,
j_s\not=0\}=J\endSb(-1)^{|J|}c_{j_1,\ldots,j_k}
=(-1)^{|J|}\[\prod_{j\in J}x_j\]f(x_1,\ldots,x_k),
\endalign$$
where in the last step we note that if $\{s:\, j_s\not=0\}=J$
and $j_s>1$ for some $s$ then $j_1+\cs+j_k>|J|\gs\deg f$ and hence
$c_{j_1,\ldots,j_k}=0$. This concludes the proof. \qed

\Remark\ 3.1. Let $f(x_1,\ldots,x_k)\in R[x_1,\ldots,x_k]$ where
$R$ is a ring with identity. It is easy to verify that
 for any $l_1,\ldots,l_k\in\N$ we have
$$\[\prod_{i=1}^k\prod_{j=1}^{l_i}x_{ij}\]
f\(\sum_{j=1}^{l_1}x_{1j},\,\ldots,\sum_{j=1}^{l_k}x_{kj}\)
=l_1!\cs l_k![x_1^{l_1}\cdots x_k^{l_k}]f(x_1,\ldots,x_k).$$
Thus, by Theorem 3.1,
$l_1!\cs l_k![x_1^{l_1}\cdots x_k^{l_k}]f(x_1,\ldots,x_k)$
is computable in terms of values of $f$ provided that $\deg f\ls l_1+\cdots+l_k$.

\proclaim{Corollary 3.1 {(\rm Escott's identity)}} Let $R$ be a ring with
identity. Given $c_1,\ldots,c_k\in R$ we have
$$\sum_{I\se[1,k]}(-1)^{|I|}\(\sum_{s\in I}c_s\)^n=0\quad\t{for every}\
n=0,1,\ldots,k-1.\tag 3.2$$
\endproclaim
\Proof. Let $n\in[0,k-1]$ and
$f(x_1,\ldots,x_k)=(\sum_{s=1}^kc_sx_s)^n$. By Theorem 3.1,
$$\sum_{I\se[1,k]}(-1)^{k-|I|}f([\![1\in I]\!],\ldots,[\![k\in I]\!])
=[x_1\cdots x_k]f(x_1,\ldots,x_k)=0.$$
This yields the desired result. \qed

\Remark\ 3.2. Escott discovered (3.2) in the case
$c_1,\ldots,c_k\in\C$ (where $\C$ is the complex field).
Maltby [M] proved Corollary 3.1
in the case where $R$ is commutative.

\proclaim{Corollary 3.2 {\rm ([Ro, Lemma 2.2])}} Let $F$ be a field,
 and let $V$ be the family of all functions from $\{0,1\}^k$ to $F$. Then
 those functions $\chi_I\in V\ (I\se[1,k])$
 given by $\chi_I(x_1,\ldots,x_k)=\prod_{s\in I}x_s$
 form a basis of the linear space $V$ over $F$.
 \endproclaim
\Proof. For $f\in V$ and $x_1,\ldots,x_k\in\{0,1\}$, clearly
$$f(x_1,\ldots,x_k)=\sum_{\da_1,\ldots,\da_k\in\{0,1\}}
f(\da_1,\ldots,\da_k)\prod_{s=1}^k[\![x_s=\da_s]\!].$$
So the dimension of $V$ does not exceed $2^k$.

 Suppose that $f=\sum_{I\se[1,k]}c_I\chi_I=0$ where $c_I\in F$.
 If $J\se[1,k]$, $c_J\not=0$ and $\deg f(x_1,\ldots,x_k)=|J|$, then
 $$c_J=\sum_{I\se J}(-1)^{|J|-|I|}f([\![1\in I]\!],\ldots,[\![k\in I]\!])=0$$
 by Theorem 3.1.
 Therefore those $\chi_I$ with $I\se[1,k]$ are linearly independent
 over $F$. We are done. \qed

 \Remark\ 3.3. Corollary 3.2 plays an important role
 in R\'onyai's study of the Kemnitz conjecture (cf. [Ro]).
  \medskip

 We mention that the Combinatorial Nullstellensatz (as stated in Section 1)
 in the important case $l_1,\ldots,l_k\in\{0,1\}$
 also follows from Theorem 3.1.
 Let $b_1\in X_1,\ldots,b_k\in X_k$ and $c_j\in X_j\sm\{b_j\}$ for $j\in
 J=\{1\ls s\ls k:\, l_s=1\}$. Set
 $$\bar f(x_1,\ldots,x_k)=f(b_1+(c_1-b_1)x_1,\ldots,b_k+(c_k-b_k)x_k)$$
 where $c_s=b_s$ for $s\in[1,k]\sm J$.
 Then $|J|=\deg f\gs\deg \bar f$ and
 $$\[\prod_{j\in J}x_j\]\bar f(x_1,\ldots,x_k)=\prod_{j\in J}(c_j-b_j)
 \times\[\prod_{j\in J}x_j\]f(x_1,\ldots,x_k)\not=0.$$
 By Theorem 3.1, for some $I\se J$ we have
 $\bar f([\![1\in I]\!],\ldots,[\![k\in I]\!])\not=0$ and hence
 $f(a_1,\ldots,a_k)\not=0$ where $a_s=b_s+(c_s-b_s)[\![s\in I]\!]\in X_s$
 for $s\in[1,k]$.

 \proclaim{Lemma 3.1 {\rm (Sun [S09])}} Let $p$ be a prime, and let $h\in\N$ and $a\in\Z$.
 Then we have the following congruence
$$\bi{a-1}{p^h-1}\eq[\![p^h\mid a]\!]\ \ (\mo\ p).\tag3.3$$
\endproclaim

Our next theorem is related to zero-sum problems
on a general abelian $p$-group $\Z_{p^{h_1}}\oplus\cdots\oplus\Z_{p^{h_l}}$.

 \proclaim{Theorem 3.2} Let
 $k,h_1,\ldots,h_l\in\Z^+$ and $k\gs\sum_{t=1}^l(p^{h_t}-1)$ where $p$
 is a prime. Let $c_{st},c_t\in\Z$ for all $s\in[1,k]$ and $t\in[1,l]$.
 Then
 $$\aligned&\sum\Sb I\se[1,k]\\p^{h_t}\mid \sum_{s\in I}c_{st}-c_t
 \ \t{for}\ t\in[1,l]\endSb(-1)^{|I|}
 \\\eq&\sum\Sb I_1\cup\cs\cup I_l=[1,k]\\|I_t|=p^{h_t}-1\ \t{for}\
 t\in[1,l]\endSb\prod_{t=1}^l\prod_{s\in I_t}c_{st}
 \ \ \ (\mo\ p).\endaligned\tag3.4$$
 \endproclaim
 \Proof. Set
 $$f(x_1,\ldots,x_k)=\prod_{t=1}^l\bi{\sum_{s=1}^kc_{st}x_s-c_t-1}{p^{h_t}-1}.$$
 Then $\deg f\ls\sum_{t=1}^l(p^{h_t}-1)\ls k$.
 Whether $n=k-\sum_{t=1}^l(p^{h_t}-1)$ is zero or not,
 $[x_1\cdots x_k]f(x_1,\ldots,x_k)$
 always coincides with
 $$[x_1\cdots x_k]\prod_{t=1}^l\f{(\sum_{s=1}^kc_{st}x_s)^{p^{h_t}-1}}{(p^{h_t}-1)!}
 =\sum\Sb I_1\cup\cs\cup I_l=[1,k]\\|I_t|=p^{h_t}-1\ \t{for}\
 t\in[1,l]\endSb\prod_{t=1}^l\prod_{s\in I_t}c_{st}.$$
 On the other hand, by Theorem 3.1 and Lemma 3.1 we have
 $$\align &(-1)^n[x_1\cdots x_k]f(x_1,\ldots,x_k)=[x_1\cdots x_k]f(x_1,\ldots,x_k)
 \\=&\sum_{I\se[1,k]}(-1)^{k-|I|}f([\![1\in I]\!],\ldots,[\![k\in I]\!])
 \\\eq&(-1)^n\sum_{I\se[1,k]}(-1)^{|I|}\prod_{t=1}^l\[\!\!\[p^{h_t}\ \bigg |
 \ \sum_{s\in I}c_{st}-c_t\]\!\!\]\ (\mo\ p).
 \endalign$$
 (Note that $(-1)^{k-n}\eq 1\ (\mo\ p)$.) Therefore (3.4) holds. \qed

 \Remark\ 3.4. In the case $k>\sum_{t=1}^l(p^{h_t}-1)$, Theorem 3.2 yields
 a theorem of Olson [O] on Davenport constants of abelian $p$-groups because the right hand side of the congruence (3.4) vanishes.
 In the same spirit, we can easily prove Theorem 2
 of Baker and Schmidt [BS] whose original proof is very deep and complicated.

 \proclaim{Corollary 3.3} Let $p$ be a prime and let $h\in\Z^+$.

 {\rm (i)} If $c,c_1,\ldots,c_{p^h-1}\in\Z$, then
 $$\sum\Sb I\se[1,p^h-1]\\p^h\mid\sum_{s\in I}c_s-c\endSb
 (-1)^{|I|}\eq c_1\cdots c_{p^h-1}\ \ (\mo\ p).\tag3.5$$

 {\rm (ii)} For $c,c_1,\ldots,c_{2p^h-2}\in\Z$ we have
 $$\aligned&\bg|\bg\{I\se[1,2p^h-2]:\, |I|=p^h-1\ \t{and}
 \ p^h\ \bigg|\ \sum_{s\in I}c_s-c\bg\}\bg|
 \\\qquad\quad&\eq[x^{p^h-1}]\prod_{s=1}^{2p^h-2}(x-c_s)\ \ (\mo\ p).
 \endaligned\tag3.6$$
 \endproclaim
 \Proof. (i) Simply apply Theorem 3.2 with $l=1$.

 (ii) In view of Theorem 3.2 in the case $l=2$,
 $$\sum\Sb I\se[1,2p^h-2]\\p^h\mid\sum_{s\in I}c_s-c
 \\p^h\mid\sum_{s\in I}1+1\endSb(-1)^{|I|}
 \eq\sum\Sb I\se[1,2p^h-2]\\|I|=p^h-1\endSb\prod_{s\in I}c_s
 \times\prod_{s\not\in I}1\ \ (\mo\ p).$$
  This is equivalent to (3.6) and we are done. \qed

 Let $q>1$ be a power of a prime $p$, and let $c_1,\ldots,c_{4q-2}\in\Z_q^2$.
 Using Lemma 3.1 and Theorem 3.1 we can prove that
 $$\align&\qquad\bg|\bg\{I\se[1,4q-2]:\, |I|=q\ \t{and}\ \sum_{s\in I}c_s=0\bg\}\bg|
 \\\eq&\bg|\bg\{I\se[1,4q-2]:\, |I|=3q\ \t{and}\ \sum_{s\in I}c_s=0\bg\}\bg|+2\ \ (\mo\ p).
 \endalign$$
This is helpful to understand the full proof of the Kemnitz conjecture given by Reiher [Re].
\medskip

 \heading{4. Proof of Theorem 2.1}\endheading

 In this section we fix a finite system (1.2) of residue classes, and set
 $I_z=\{1\ls s\ls k:\, z\in a_s(n_s)\}$ for $z\in\Z$. We first extend [S09, Lemma 4.1]
 to any field containing an element of $($multiplicative$)$ order $N_A$, where $N_A$
 is the least common multiple of the moduli $n_1,\ldots,n_k$ in $(1.2)$.

 \proclaim{Lemma 4.1} Let $A$ be as in $(1.2)$ and let $m_1,\ldots,m_k\in\Z$.
Let $F$ be a field containing an element $\zeta$ of $($multiplicative$)$ order $N_A$,
and let $f(x_1,\ldots,x_k)$ be a polynomial over $F$ with $\deg f\ls m(A)$. If
$[\prod_{s\in I_z}x_s]f(x_1,\ldots,x_k)=0$
for all $z\in\Z$, then we have $\psi(\theta)=0$ for any
$0\ls \theta<1$, where
$$\psi(\theta):=\sum\Sb I\se[1,k]\\\{\sum_{s\in I}m_s/n_s\}=\theta\endSb
(-1)^{|I|}f([\![1\in I]\!],\ldots,[\![k\in I]\!])\zeta^{N_A\sum_{s\in I}a_sm_s/n_s}.$$
The converse holds when $m_1,\ldots,m_k$ are relatively prime to $n_1,\ldots,n_k$
respectively.
\endproclaim
\Proof. Let $z\in\Z$ and $J\se[1,k]$. Clearly
$$\align&[\![J\supseteq I_z]\!]\prod_{s=1}^k\l([\![s\not\in J]\!]-\zeta^{N_A(a_s-z)m_s/n_s}\r)
\\=&\sum_{I\se[1,k]}\prod^k\Sb s=1\\s\not\in I\endSb
[\![s\not\in J]\!]\times (-1)^{|I|}
\zeta^{N_A\sum_{s\in I}a_sm_s/n_s}\zeta^{-zN_A\sum_{s\in I}m_s/n_s}
\\=&\sum_{\theta\in S}\zeta^{-zN_A\theta}
\sum\Sb J\se I\se[1,k]\\\{\sum_{s\in I}m_s/n_s\}=\theta\endSb
(-1)^{|I|}\zeta^{N_A\sum_{s\in I}a_sm_s/n_s}
\endalign$$
where
$$S=\bg\{\bg\{\sum_{s\in I}\f{m_s}{n_s}\bg\}:\, I\se[1,k]\bg\}.\tag4.1$$

Write
$f(x_1,\ldots,x_k)=\sum_{j_1,\ldots,j_k\gs0}c_{j_1,\ldots,j_k}x_1^{j_1}\cdots x_k^{j_k}$.
Obviously
$$f([\![1\in I]\!],\ldots,[\![k\in I]\!])=\sum\Sb j_1,\ldots,j_k\gs0
\\\{1\ls s\ls k:\, j_s\not=0\}\se I\endSb c_{j_1,\ldots,j_k}
\quad\t{for all}\ I\se[1,k].$$
If $c_{j_1,\ldots,j_k}\not=0$ and $J=\{1\ls s\ls k:\, j_s\not=0\}\supseteq I_z$, then
$\deg f\gs|J|\gs|I_z|=w_A(z)\gs\deg f$;
hence $w_A(z)=\deg f$, $J=I_z$ and $j_s=1$ for $s\in J$.

In view of the above, for any $z\in\Z$ the sum $\sum_{\theta\in S}\zeta^{-zN_A\theta}\psi(\theta)$
coincides with
$$\align&\sum_{\theta\in S}\zeta^{-zN_A\theta}\sum\Sb I\se[1,k]\\\{\sum_{s\in I}m_s/n_s\}=\theta\endSb
(-1)^{|I|}\sum\Sb j_1,\ldots,j_k\gs0
\\\{s:\, j_s\not=0\}\se I\endSb c_{j_1,\ldots,j_k}\zeta^{N_A\sum_{s\in I}a_sm_s/n_s}
\\=&\sum_{j_1,\ldots,j_k\gs0}c_{j_1,\ldots,j_k}
\sum_{\theta\in S}\zeta^{-zN_A\theta}
\sum\Sb \{s:\, j_s\not=0\}\se I\se[1,k]\\\{\sum_{s\in I}m_s/n_s\}=\theta\endSb
(-1)^{|I|}\zeta^{N_A\sum_{s\in I}a_sm_s/n_s}
\\=&\sum\Sb j_1,\ldots,j_k\gs0\\ J=\{s:\,j_s\not=0\}\sp I_z\endSb c_{j_1,\ldots,j_k}
\prod_{s=1}^k\l([\![s\not\in J]\!]-\zeta^{N_A(a_s-z)m_s/n_s}\r)
\\=&c(I_z)\prod_{s=1}^k\l([s\not\in I_z]-\zeta^{N_A(a_s-z)m_s/n_s}\r),
\endalign$$
where $c(I_z)=[\prod_{s\in I_z}x_s]f(x_1,\ldots,x_k)$. Therefore
$$\sum_{\theta\in S}\zeta^{-zN_A\theta}\psi(\theta)
=(-1)^kc(I_z)\prod^k\Sb s=1\\s\not\in I_z\endSb\l(\zeta^{N_A(a_s-z)m_s/n_s}-1\r).\tag4.2$$
When $n_1=\cs=n_k=1$, this yields the equality
$$\sum_{I\se[1,k]}(-1)^{|I|}f([\![1\in I]\!],\ldots,[\![k\in I]\!])=(-1)^k[x_1\cs x_k]f(x_1,\ldots,x_k)$$
as asserted by Theorem 3.1.

Observe that $c(I_z)=0$
if $\psi(\theta)=0$ for all $0\ls\theta<1$
and each $m_s$ is relatively prime to $n_s$.

Suppose that $c(I_z)=0$ for all $z\in\Z$. Then
$\sum_{\theta\in S}\zeta^{-nN_A\theta}\psi(\theta)=0$ for all
$n\in[0,|S|-1]$. As the Vandermonde-type determinant
$$\det[(\zeta^{-N_A\theta})^n]_{n\in[0,|S|-1],\ \theta\in S}$$
is nonzero,
we have $\psi(\theta)=0$ for all $\theta\in S$.
If $0\ls\theta<1$ and $\theta\not\in S$, then
$\psi(\theta)=0$ holds trivially.

 In view of the above, we have completed the proof of Lemma 4.1. \qed

\proclaim{Lemma 4.2} Let $F$ be a field of characteristic $p$,
and let $n$ be a positive integer. Then
$p\nmid n$ if and only if there is an extension field of $F$
containing an element of $($multiplicative$)$ order $n$.
\endproclaim
\Proof. (i) Suppose that $p\mid n$ and $E/F$ is a field extension.
If $\zeta\in E$ and $\zeta^n=1$, then
$(\zeta^{n/p}-1)^p=(\zeta^{n/p})^p-1=\zeta^n-1=0$
and hence $\zeta^{n/p}-1=0$. So $E$ contains no element of order
$n$.

(ii) Now assume that $p\nmid n$. Let $E$ be the splitting field of
the polynomial $f(x)=x^{n}-1$ over $F$.
Then $G=\{\zeta\in E:\, \zeta^{n}=1\}$
is a finite subgroup of the multiplicative group $E^*=E\sm\{0\}$,
therefore it is cyclic by field theory.
Since $p\nmid n$, $f'(\zeta)=n\zeta^{n-1}\not=0$
for any $\zeta\in G$. So the equation $f(x)=0$ has no repeated roots in $E$
and hence $|G|=n$. Any generator of the cyclic group $G$
has order $n$.

Combining the above we obtain the desired result. \qed

\medskip
\noindent{\tt Proof of Theorem 2.1}.
For convenience we set $h=|J|-\deg P$, $A=\{a_s(n_s)\}_{s=1}^k$,
$J^*=\{1\ls s\ls k:\, a_0\in a_s(n_s)\}$
and $J'=J^*\sm J$.

Let $d_1,\ldots,d_h$ be any elements of $F$
and define
$$\align f(x_1,\ldots,x_k)=&P(b_1+(c_1-b_1)x_1,\ldots,b_k+(c_k-b_k)x_k)
\\&\times\prod_{j=1}^h\(\sum_{s=1}^k(b_s+(c_s-b_s)x_s)-d_j\)\times\prod_{s\in J'}(x_s-1).
\endalign$$
Then $\deg f\ls\deg P+|J'|+h=|J^*|=m(A)$. As
$[\prod_{s\in J^*}x_s]f(x_1,\ldots,x_k)$ equals
$$\prod_{j\in J}(c_j-b_j)\times\[\prod_{j\in J}x_j\]P(x_1,\ldots,x_k)
\(\sum_{s=1}^kx_s\)^h\not=0,$$
we have $\deg f=m(A)$. Recall that $m_s$ is relatively prime
to $n_s$ for each $s\in[1,k]$.
In light of Lemma 4.1,
$$\psi(\theta)=\sum\Sb I\se[1,k]\\\{\sum_{s\in I}m_s/n_s\}=\theta\endSb
(-1)^{|I|}f([\![1\in I]\!],\ldots,[\![k\in I]\!])\zeta^{N_A\sum_{s\in I}a_sm_s/n_s}\not=0$$
for some $0\ls\theta<1$,
where $\zeta$ is an element of order $N=N_A$
in an extension field of $F$ (whose existence follows from Lemma 4.2).

Let $\al=\{n_0\theta\}$ and $r\in[0,n_0-1]$.
Then $(\al+r)/{n_0}=\theta+\bar r/n_0$ where $\bar r=r-\lfloor n_0\theta\rfloor$.
As $w_A(a_0+N)=w_A(a_0)<w_{A_0}(a_0)=m(A_0)$,
we must have $a_0+N\in a_0(n_0)$ and hence $n_0\mid N$.
Note that $\deg f=m(A)<m(A_0)$.
Applying Lemma 4.1 to the system $A_0$ we find that
$$\psi(\theta)+\sum\Sb I\se[1,k]\\\{\sum_{s\in I}m_s/n_s+(-\bar r)/{n_0}\}
=\theta\endSb(-1)^{|I|+1}f([\![1\in I]\!],\ldots,[\![k\in I]\!])
\zeta^{N\beta(I)}=0$$
where $\beta(I)=\sum_{s\in I}a_sm_s/n_s+a_0(-\bar r)/n_0$.
It follows that
$$\psi\l(\f{\al+r}{n_0}\r)=\psi\l(\theta+\f{\bar r}{n_0}\r)
=\zeta^{Na_0\bar r/n_0}\psi(\theta)\not=0.$$
So, there exists an $I\se[1,k]$ with $\{\sum_{s\in I}m_s/n_s\}=(\al+r)/{n_0}$
such that $f([\![1\in I]\!],\ldots,[\![k\in I]\!])\not=0$.
Note that $x_s=b_s+(c_s-b_s)[\![s\in I]\!]\in X_s$ for all $s\in[1,k]$.
Also, $P(x_1,\ldots,x_k)\not=0$, $I\cap J'=\em$
and $I=\{1\ls s\ls k:\, x_s\not=b_s\}$.
Thus $S_r$ contains
$x_1+\cs+x_k$
which is different from $d_1,\ldots,d_h$.

If $|S_r|\ls h$, then we can select
$d_1,\ldots,d_h\in F$ such that $\{d_1,\ldots,d_h\}=S_r,$
 hence we get a contradiction from the above.
 Therefore $|S_r|\gs h+1$ and we are done. \qed

\enddocument